\documentclass[titlepage,twoside,12pt]{article}
\usepackage{amssymb}
\usepackage{amsfonts}
\begin{document}

{\LARGE \bf New Method for Solving Large Classes \\ \\ of Nonlinear Systems of PDEs } \\ \\

{\bf Elem\'{e}r E ~Rosinger} \\ \\
{\small \it Department of Mathematics \\ and Applied Mathematics} \\
{\small \it University of Pretoria} \\
{\small \it Pretoria} \\
{\small \it 0002 South Africa} \\
{\small \it eerosinger@hotmail.com} \\ \\

{\bf Abstract} \\

The essentials of a new method in solving very large classes of nonlinear systems of PDEs,
possibly associated with initial and/or boundary value problems, are presented. The PDEs can
be defined by continuous, not necessarily smooth expressions, and the solutions obtained cab be
assimilated with usual measurable functions, or even with Hausdorff continuous ones. The
respective result sets aside completely, and with a large nonlinear margin, the celebrated
1957 impossibility of Hans Lewy regarding the nonexistence of solution in distributions of
large classes of linear smooth coefficient PDEs.  \\ \\ \\

{\bf 1. The Class of Nonlinear Systems of PDEs Solved} \\

The nonlinear systems of PDEs considered in this paper are composed of equations of the
general form \\

(1.1)~~~ $ F ( x, U ( x ), ~.~.~.~ , D^p_x U ( x ), ~.~.~.~ ) ~=~ f ( x ),~~~
                                               x \in \Omega ~\subseteq~ \mathbb{R}^n $ \\

were the domains $\Omega$ can be any open, not necessarily bounded subsets of $\mathbb{R}^n$,
while $p \in \mathbb{N}^n,~ | p | \leq m$, with the orders $m \in \mathbb{N}$ of the PDEs
arbitrary given. \\

The functions $F$ which define the left hand terms are only assumed to be {\it jointly
continuous} in all of their arguments. The right hand terms $f$ are also required to be {\it
continuous} only. \\
However, with minimal modifications of the method, both $F$ and $f$ can have certain {\it
discontinuities} as well, [3]. \\

As it turns out, regardless of the above generality of the nonlinear systems of PDEs
considered, and of possibly associated initial and/or boundary value problems, one can always
find for them solutions $U$ defined on the {\it whole} of the respective domains $\Omega$.
These solutions $U$ have the {\it blanket, type independent}, or {\it universal regularity}
property that they can be assimilated with usual {\it measurable}, or even {\it Hausdorff
continuous functions}, [1,4-6]. \\

Thus when solving systems of nonlinear PDEs of the generality of those in (1.1), one can {\it
dispense with} the various customary spaces of distributions, hyperfunctions, generalized
functions, Sobolev spaces, and so on. Instead one can stay within the realms of {\it usual
functions}. However, functional analytic methods can possibly be used in order to obtain
further regularity or other desirable properties of such solutions. \\ \\

{\bf 2. A Natural Smooth Framework for Nonlinear \\
        \hspace*{0.5cm} Systems of PDEs} \\

Let us now associate with a nonlinear PDE in (1.1) the corresponding nonlinear partial
differential operator defined by the left hand side, namely \\

(2.1)~~~ $ T ( x, D ) U ( x ) ~=~ F ( x, U ( x ), ~.~.~.~ , D^p_x U ( x ), ~.~.~.~ ),~~~
                                        x \in \Omega $ \\

{\it Two} facts about the nonlinear PDEs in (1.1) and the corresponding nonlinear partial
differential operators $T ( x, D )$ in (2.1) are important and immediate :

\begin{itemize}

\item The operators $T ( x, D )$ can {\it naturally} be seen as acting in the {\it classical}
or smooth context, namely

\end{itemize}

(2.2)~~~ $ T ( x, D ) ~:~ {\cal C}^m ( \Omega ) \ni U ~~\longmapsto~~ T ( x, D ) U \in
                                          {\cal C}^0 ( \Omega ) $ \\

while, unfortunately on the other hand :

\begin{itemize}

\item The mappings in this natural classical context (2.2) are typically {\it not} surjective
even in the case of linear $T ( x, D )$, and they are even less so in the general nonlinear
case of (1.1).

\end{itemize}

In other words, linear or nonlinear PDEs in (1.1) typically {\it cannot} be expected to have
{\it classical} solutions $U \in {\cal C}^m ( \Omega )$, for arbitrary continuous right hand
terms $f \in {\cal C}^0 ( \Omega )$, as illustrated by a variety of well known examples, some
of them rather simple ones, see [3, chap. 6]. \\
Furthermore, it can often happen that nonclassical solutions do have a major applicative
interest, thus they have to be sought out {\it beyond} the confines of the classical framework
in (2.2). \\
This is, therefore, how we are led to the {\it necessity} to consider {\it generalized
solutions} $U$ for PDEs like those in (1.1), that is, solutions $U \notin {\cal C}^m ( \Omega
)$, which therefore are no longer classical. This means that the natural classical mappings
(2.2) must in certain suitable ways be {\it extended} to {\it commutative diagrams}

\begin{math}
\setlength{\unitlength}{0.2cm}
\thicklines
\begin{picture}(60,20)

\put(10,15){${\cal C}^m ( \Omega )$}
\put(27,17){$T ( x, D )$}
\put(18,15.5){\vector(1,0){26.5}}
\put(47,15){${\cal C}^0 ( \Omega )$}
\put(0,8){$(2.3)$}
\put(12,13){\vector(0,-1){8}}
\put(13.5,8){$\subseteq$}
\put(49,13){\vector(0,-1){8}}
\put(45.7,8){$\subseteq$}
\put(10,2){$~~X$}
\put(16,2.6){\vector(1,0){29.5}}
\put(47,2){$~~Y$}
\put(29,-0.5){$\overline{T}$}

\end{picture}
\end{math}

\smallskip
with the generalized solutions now being found as \\

(2.4) ~~~$ U \in X \setminus {\cal C}^m ( \Omega ) $ \\

instead of the classical ones $U \in {\cal C}^m ( \Omega )$ which may easily fail to exist. A
further important point is that one expects to reestablish certain kind of {\it surjectivity}
type properties typically missing in (2.2), at least such as for instance \\

(2.5)~~~ $ {\cal C}^0 ( \Omega ) ~\subseteq~ \overline{T} ( X ) $ \\

As it turns out, when constructing extensions of (2.2) given by commutative diagrams (2.3), we
shall be interested in the following somewhat larger spaces of piecewise smooth functions. For
any integer $0 \leq l \leq \infty$, we define \\

(2.6)~~~ $ {\cal C}^l_{nd} ( \Omega ) ~=~ \left \{~ u : \Omega \rightarrow \mathbb{R}~~
                    \begin{array}{|l}
                        ~\exists~ \Gamma \subset \Omega ~\mbox{closed, nowhere dense}~ : \\
                        ~~~~ u \in {\cal C}^l ( \Omega \setminus \Gamma )
                     \end{array} ~\right \} $ \\

and as an immediate strengthening of (2.2), we obviously obtain \\

(2.7)~~~ $ T ( x, D )~ {\cal C}^m_{nd} ( \Omega ) ~\subseteq ~
                                         {\cal C}^0_{nd} ( \Omega ) $ \\

The solution of the nonlinear PDEs in (1.1) through the order completion method will come from
the construction of specific instances of the {\it commutative diagrams} (2.3), given by \\

\begin{math}
\setlength{\unitlength}{0.2cm}
\thicklines
\begin{picture}(60,30)

\put(10,25){${\cal C}^m_{nd} ( \Omega )$}
\put(27,27){$T ( x, D )$}
\put(18,25.5){\vector(1,0){26.5}}
\put(47,25){${\cal C}^0_{nd} ( \Omega )$}
\put(0,14){$(2.8)$}
\put(12,23){\vector(0,-1){6}}
\put(49,23){\vector(0,-1){6}}
\put(10,14){${\cal M}^m_T ( \Omega )$}
\put(19,14.5){\vector(1,0){25.5}}
\put(47,14){${\cal M}^0 ( \Omega )$}
\put(29,16){$T$}
\put(12,12){\vector(0,-1){6}}
\put(49,12){\vector(0,-1){6}}
\put(10,3){${\cal M}^m_T ( \Omega )^{\#}$}
\put(47,3){${\cal M}^0 ( \Omega )^{\#}$}
\put(26.5,5){$\mbox{bijective}$}
\put(19,3.5){\vector(1,0){25.5}}
\put(30,0.5){$T^{\#}$}

\end{picture}
\end{math} \\

where the operation $(~~)^{\#}$ means the {\it order completion}, [3], of the respective
spaces, as well as the extension to such order completions of the respective mappings, see
(5.7) below. It follows that in terms of (2.3), we have \\

$ X ~=~ {\cal M}^m_T ( \Omega )^{\#},~~~ Y ~=~ {\cal M}^0 ( \Omega )^{\#},~~~
                                                      \overline{T} ~=~ T^{\#} $ \\

thus we shall obtain for the nonlinear PDEs in (1.1) generalized solutions \\

(2.9) ~~~$ U \in {\cal M}^m_T ( \Omega )^{\#} $ \\

Furthermore, instead of the {\it surjectivity} condition (2.5), we shall at least have the
following stronger one \\

(2.10) ~~~$ {\cal C}^0_{nd} ( \Omega ) ~\subseteq~
                     T^{\#} ( {\cal M}^m_T ( \Omega )^{\#} ) $ \\

So far about the main ideas related to the {\it existence} of solutions of general nonlinear
PDEs of the form (1.1). Further details can be found in [3,1,4-6]. \\

As for the {\it regularity} of such solutions, we recall that, as shown in [1], one has the
inclusions \\

(2.11) ~~~$ {\cal M}^0 ( \Omega )^{\#} ~\subseteq~ Mes\, ( \Omega ) $ \\

where $Mes\, ( \Omega )$ denotes the set of Lebesgue measurable functions on $\Omega$. In this
way, in view of (2.8) and (2.9), one can assimilate the generalized solutions $U$ of the
nonlinear PDEs in (1.1) with usual measurable functions in $Mes\, ( \Omega )$. \\
Recently, however, based on results in [1,4-6], it was shown that instead of (2.11), one has
the much {\it stronger regularity} property \\

(2.12) ~~~$ {\cal M}^0 ( \Omega )^{\#} ~\subseteq~ \mathbb{H}\, ( \Omega ) $ \\

where $\mathbb{H}\, ( \Omega )$ denotes the set of Hausdorff continuous functions on $\Omega$.
Consequently, now one can significantly improve on the earlier regularity result, as one can
assimilate the generalized solutions $U$ of the nonlinear PDEs in (1.1) with usual functions
in $\mathbb{H}\, ( \Omega )$. \\

Regarding {\it systems} of nonlinear PDEs such as in (1.1), with possibly associated initial
and/or boundary value problems, it was shown in [3, chap. 8] the way they can be dealt with
the above order completion method. \\
In this respect, a {\it surprising} advantage of the order completion method is the ease, when
compared with the usual functional analytic approaches, in dealing with initial and/or
boundary value problems. \\ \\

{\bf 3. Solution Method by Order Completion} \\

As first introduced and developed in [3], the mentioned {\it existence} and {\it blanket, type
independent}, or {\it universal regularity} property of solutions for nonlinear systems of
PDEs of the above generality in (1.1), together with possibly associated initial and/or
boundary value problems, is based on a rather simple, however, quite general {\it order
completion} of suitable spaces of piecewise smooth functions. \\

Here one can mention that, in fact, this order completion based solution method reaches {\it
far beyond} the solution of PDEs, and in fact it can be applied to the solution of general
equations \\

(3.1)~~~ $ T ( A ) ~=~ F $ \\

where \\

(3.2)~~~ $ T : X ~\longrightarrow~ Y $ \\

is any mapping, $X$ is any nonvoid set, while $( Y, \leq )$ is a partially ordered set, or in
short, poset, while $F \in Y$ is given, and $A \in X$ is the sought after solution. \\

In general, for a given $F \in Y$, there may not exist any solution $A \in X$ for the equation
(3.2), unless the setup (3.1), (3.2) is {\it extended} in suitable ways. \\
The usual such extensions assume suitable topologies on $X$ and $Y$, and certain continuity
properties for the mapping $T$ in (3.2). \\
As shown in [3], and also seen in the sequel, for the same purpose of solving the equations
(3.1) in a suitably extended setup, one can successfully use the {\it order completion} of
the spaces $X$ and $Y$. \\

A {\it major advantages} of the order completion approach is that such a method

\begin{itemize}

\item does {\it no longer} differentiate between linear and nonlinear operators $T$ in (3.2),
in case the spaces $X$ and $Y$ may happen to have a linear vector space structure, [3].

\end{itemize}

Two further convenient features of the order completion method in solving general equations
(3.1) are the following :

\begin{itemize}

\item one obtains necessary and sufficient conditions for the existence of solutions,

\item one obtains explicit expressions of the solutions, whenever they exist.

\end{itemize}

The order completion method for solving general equations (3.1) presented here is but one of
the possible such approaches. This method is based on so called "pull-back" partial orders on
the domains of the mappings $T$ in (3.1), (3.2). However, the order completion method for
solving general equations (3.1) is {\it not} limited to "pull-back" partial orders. Further
possible developments in this regard of the order completion method, and {\it no longer} based
on "pull-back" partial orders will be indicated. \\

Here it should be mentioned that usual functional analytic methods for solving linear or
nonlinear PDEs often make use of "pull-back" topologies on the domains of the respective
partial differential operators as illustrated next in section 4, see for further details in
[3, chap. 12]. \\ \\

{\bf 4. Use of "Pull-Back" in Functional Analytic Solution \\
        \hspace*{0.5cm} Methods of PDEs} \\

We present here one of the classical examples where a "pull-back" topology on the domain of a
partial differential operator is used in order to solve the corresponding PDE. Let us consider
on a bounded Euclidean domain $\Omega$, which has a smooth boundary $\partial~ \Omega$, the
following familiar linear boundary value problem, usually called the Poisson Problem \\

(4.1) $~~~ \begin{array}{l}
                    \Delta~ U(x) ~=~ f(x),~~~ x \in \Omega \\ \\
                     U ~=~ 0 ~~~\mbox{on}~~ \partial~ \Omega
             \end{array} $ \\

As is well known, for every given $f \in C^\infty ( \overline{\Omega} )$, where
$\overline{\Omega}$ denotes the closure of $\Omega$, this problem has a unique solution $U$
in the space \\

(4.2) $~~~ X ~=~ \left \{~~ v \in C^\infty (\overline{\Omega}) ~~|~~ v ~=~ 0 ~~\mbox{on}~~
                                                            \partial~ \Omega ~~\right \} $ \\

It follows that the mapping \\

(4.3) $~~~ X \ni v ~\longmapsto~ | |~ \Delta v ~| |_{~L^2 (~\Omega~)} $ \\

defines a norm on the vector space $X$. Now let \\

(4.4) $~~~ Y ~=~ C^\infty (\overline{\Omega}) $ \\

be endowed with the topology induced by $L^2 (\Omega)$. Then in view of (4.1) - (4.4), the
mapping  \\

(4.5) $~~~ \Delta : X ~\rightarrow~ Y $ \\

is a uniform continuous linear bijection. Therefore, it can be extended in a unique manner to
an isomorphism of Banach spaces \\

(4.6) $~~~ \Delta : \overline{X} ~\rightarrow~ \overline{Y} ~=~  L^2 (\Omega) $ \\

In this way one has the classical existence and uniqueness result \\

(4.7) $~~~ \begin{array}{l}
                  \forall~~~ f \in L^2 (\Omega)~~ : \\ \\
                  \exists~!~~ U \in \overline{X} ~~: \\ \\
                   ~~~~~~ \Delta U = f
              \end{array} $ \\

The power and simplicity - based on linearity and topological completion of uniform spaces -
of the above classical existence and uniqueness result is obvious. This power is illustrated
by the fact that the set $\overline{Y} =  L^2 (\Omega)$ in which the right hand terms $f$ in
(4.1) can now be chosen is much {\it larger} than the original $Y = C^\infty
(\overline{\Omega})$.
Furthermore, the existence and uniqueness result in (4.7) does not need the a priori knowledge
of the structure of the elements $U \in \overline{X}$, that is, of the respective generalized
solutions. This structure which gives the regularity properties of such solutions can be
obtained by a further detailed study of the respective differential operators defining the
PDEs under consideration, in this case, the Laplacean $\Delta$. And in the above specific
instance we obtain \\

(4.8) $~~~ \overline{X} ~=~ H^2 (\Omega) \cap H^1_0 (\Omega) $ \\

As seen above, typically for the functional analytic methods, the generalized solutions are
obtained in topological completions of vector spaces of usual functions. And such completions,
like for instance the various Sobolev spaces, are defined by certain linear partial
differential operators which may happen to {\it depend} on the PDEs under consideration. \\

In the above example, for instance, the topology on the space $X$ obviously {\it depends} on
the specific PDE in (4.1). Thus the topological completion $\overline{X}$ in which the
generalized solutions $U$ are found according to (4.7), does again {\it depend} on the
respective PDE. \\ \\

{\bf 5. Pull-Back Partial Order} \\

Without loss of generality, [3], we shall assume that all the posets considered are without a
minimum or maximum element. Various notions and results related to partial orders which are
used in the sequel are presented in the Appendix. \\

Given an equation (3.1), (3.2), we define on $X$ the equivalence relation $\approx_T$ by \\

(5.1)~~~ $ u \approx_T y ~~~\Longleftrightarrow~~~ T ( u ) ~=~ T ( v ) $ \\

for $u, v \in X$. In this way, by considering the quotient space \\

(5.2)~~~ $ X_T ~=~ X / \approx_T $ \\

we obtain the {\it injective} mapping \\

(5.3)~~~ $ T_\approx : X_T ~\longrightarrow~ Y $ \\

defined by \\

(5.4)~~~ $ X_T \ni U \longmapsto T ( u ) \in Y $ \\

where $u \in U$, that is, $U$ is the $\approx_T$ equivalence class of $u$ in $X_T$, while
$T ( u )$ is defined by (3.2). \\

At that stage, we can define a partial order $\leq_T$ on $X_T$ as being the {\it pull-back} by
the mapping $T_\approx$ in (5.3) of the given partial order $\leq$ on $Y$, namely \\

(5.5)~~~ $ U ~\leq_T~ V ~~~\Longleftrightarrow~~~ T_\approx ( U ) ~\leq~ T_\approx ( V ) $ \\

for $U, V \in X_T$. The effect of the above construction is that we obtain the {\it order
isomorphic embedding}, or in short OIE \\

(5.6)~~~ $ T_\approx : X_T ~\longrightarrow~ Y $ \\

As mentioned, without loss of generality we shall assume that the poset $( X^{\#}_T, \leq_T )$
has no minimum or maximum. \\

And now, we consider the {\it order completions} $X^{\#}_T$ and $Y^{\#}$ of $X_T$, and
respectively, $Y$. \\
For simplicity, we shall denote by $\leq$ the partial orders both on $X^{\#}_T$ and $Y^{\#}$.
In fact, as seen in (A.20), these partial orders are the usual inclusion relations $\subseteq$
among subsets of $X$, respectively, of $Y$. \\

Then according to Proposition A.1 in the Appendix, we obtain the commutative diagram of
OIE-s \\

\begin{math}
\setlength{\unitlength}{0.2cm}
\thicklines
\begin{picture}(60,19)

\put(11,15){$X_T$}
\put(29,17){$T_\approx$}
\put(16,15.5){\vector(1,0){28.5}}
\put(47,15){$Y$}
\put(0,8){$(5.7)$}
\put(12,13){\vector(0,-1){8}}
\put(13.5,8){$\subseteq$}
\put(47.5,13){\vector(0,-1){8}}
\put(44.5,8){$\subseteq$}
\put(11,2){$X^{\#}_T$}
\put(16,2.6){\vector(1,0){28.5}}
\put(47,2){$Y^{\#}$}
\put(29,-0.5){$T^{\#}$}

\end{picture}
\end{math} \\

Consequently, for $U \in X_T$ and $A \in X^{\#}_T$, we have in $Y^{\#}$ the relations \\

(5.8)~~~ $ T^{\#} ( < U\, ]\,) ~=~ < T_\approx ( U )\, ] $ \\

(5.9)~~~ $ T^{\#} ( A ) ~=~ ( T_\approx ( A ) )^{ul} ~=~
                        \sup_{\,Y^{\#}}~ \{~ < T_\approx ( U )\,] ~~|~~ U \in A ~\} $ \\ \\

{\bf 6. Reformulation} \\

Now we can reformulate the problem of solving the general equations (3.1), (3.2) as follows.
Given $F \in Y^{\#}$, find {\it necessary and sufficient} conditions for the existence of $A
\in X^{\#}_T$, such that \\

(6.1) $~~~ T^{\#} ( A ) ~=~ F $ \\ \\

{\bf 7. Solution} \\

We note that (5.7) gives the inclusions \\

(7.1) $~~~ \begin{array}{l}
               \sup_{Y^{\#}}~ \{~ T^{\#} ( U ) ~|~ U \in X^{\#}_T,~
                          T^{\#} ( U ) \subseteq F ~\} ~\subseteq~ \\ \\
               \subseteq~ T^{\#} (~ \sup_{X^{\#}_T}~ \{~ U ~|~ U \in X^{\#}_T,~
                          T^{\#} ( U ) \subseteq F ~\} ~) ~\subseteq~ \\ \\
               \subseteq~ T^{\#} (~ \inf_{X^{\#}_T}~ \{~ V ~|~ V \in X^{\#}_T,~
                          F ~\subseteq~ T^{\#} ( V ) ~\} ~) ~\subseteq~ \\ \\
               \subseteq~ \inf_{Y^{\#}}~ \{~ T^{\#} ( V ) ~|~ V \in X^{\#}_T,~
                          F ~\subseteq~ T^{\#} ( V ) ~\}
            \end{array} $ \\ \\

Indeed, the first and last inclusions follow from Lemma A.1 in the Appendix. As for the middle
inclusion in (7.1), let $U, V \in X^{\#}_T$ be such that $T^{\#} ( U ) \subseteq F
\subseteq T^{\#} ( V )$. Then $T^{\#} ( U ) \subseteq T^{\#} ( V )$, hence $U \leq V$,
since $T^{\#}$ is an OIE. It follows that \\

$~~~ \sup_{X^{\#}_T}~ \{~ U ~|~ U \in X^{\#}_T,~ T^{\#} ( U ) \subseteq F ~\} ~\leq~ \\ \\
 ~~~~~~\leq~ \inf_{X^{\#}_T}~ \{~ V ~|~ V \in X^{\#}_T,~ F ~\subseteq~ T^{\#} ( V ) ~\} $ \\

and the proof of (7.1) is completed. \\

We note further that the above inequality $T^{\#} ( U ) \subseteq F \subseteq T^{\#} ( V )$
also implies \\

(7.2) $~~~ \begin{array}{l}
                    \sup_{Y^{\#}}~ \{~ T^{\#} ( U ) ~|~ U \in X^{\#}_T,~ T^{\#} ( U )
                             \subseteq F ~\} ~\subseteq~ F ~\subseteq~ \\ \\
           ~~~~~~~~~~~~\subseteq~ \inf_{Y^{\#}}~ \{~ T^{\#} ( V ) ~|~ V \in X^{\#}_T,~ F
                             ~\subseteq~ T^{\#} ( V ) ~\}
            \end{array} $ \\

Furthermore \\

(7.3) $~~~ \{~ U \in X^{\#}_T ~|~ T^{\#} ( U ) ~\subseteq~ F ~\} ~\neq~ \phi $ \\

since (A.9) gives $U = \phi \in X^{\#}_T$, hence in view of (A.4) - (A.6) and (5.9) we have
$T^{\#} ( U ) = ( T ( \phi ) )^{ul} = \phi^{ul} = ( \phi^u )^l = ( Y^{\#} )^l = \phi
\subseteq F$. \\

Returning now to the problem (6.1), we note that it is not trivial. Indeed, the OIE in (5.7),
namely  \\

$~~~ T^{\#} : X^{\#}_T ~\longrightarrow~ Y^{\#} $ \\

need {\it not} be surjective. Further $T^{\#}$ need {\it not} preserve infima or suprema. \\

However, in the next theorem we can obtain the following two general results :

\begin{itemize}

\item a necessary and sufficient condition for the solvability of (6.1), and

\item the explicit expression of the solution, when it exists.

\end{itemize}

\bigskip
{\bf Theorem 7.1.} \\

Given $F \in Y^{\#}$. \\

1)~ The equation \\

(7.4) $~~~ T^{\#} ( A ) ~=~ F $ \\

has a solution $A \in X^{\#}_T$, if and only if, see (7.1) \\

(7.5) $~~~ \begin{array}{l}
                      \sup_{Y^{\#}}~ \{~ T^{\#} ( U ) ~|~ U \in X^{\#}_T,~ T^{\#} ( U )
                                   \subseteq F ~\} ~=~ \\ \\
   ~~~~~~~~~~~~~~~~~~~~~~~~=~ \inf_{Y^{\#}}~ \{~ T^{\#} ( V ) ~|~ V \in X^{\#}_T,~
                                    F ~\subseteq~ T^{\#} ( V ) ~\}
            \end{array} $ \\

2)~ This solution is unique, whenever it exists, see (5.7). \\

3)~ When it exists, the unique solution $A \in X^{\#}_T$ is given by \\

(7.6) $~~~ \begin{array}{l}
                   A ~=~ \sup_{X^{\#}_T}~ \{~ U \in X^{\#}_T ~|~ T^{\#} ( U )
                                      ~\subseteq~ F ~\} ~=~ \\ \\
            ~~~~~=~ \inf_{X^{\#}_T}~ \{~ V \in X^{\#}_T ~|~ F ~\subseteq~ T^{\#} ( V ) ~\}
             \end{array} $ \\

and, see (7.3) \\

(7.7) $~~~ \{~ U \in X^{\#}_T ~|~ T^{\#} ( U ) ~\subseteq~ F ~\},~~
            \{~ V \in X^{\#}_T ~|~ F ~\subseteq~ T^{\#} ( V ) ~\} ~\neq~ \phi $ \\

{\bf Proof} \\

From (7.4) follows that $A \in X^{\#}_T,~ T^{\#} ( A ) \subseteq F$, thus \\

$~~~ F ~=~ T^{\#} ( A ) ~\subseteq~ \sup_{Y^{\#}}~ \{~ T^{\#} ( U ) ~|~ U \in X^{\#}_T,~
                          T^{\#} ( U ) \subseteq F ~\} $ \\

Similarly we have \\

$~~~ \inf_{Y^{\#}}~ \{~ T^{\#} ( V ) ~|~ V \in X^{\#}_T,~
                          F ~\subseteq~ T^{\#} ( V ) ~\} ~\subseteq~ T^{\#} ( A ) ~=~ F $ \\

thus (7.1) collapses to the seven equalities \\

$~~~ F ~=~ T^{\#} ( A ) ~=~ \sup_{Y^{\#}}~.~.~.~ ~=~
                    T^{\#} ( \sup_{X^{\#}_T}~.~.~.~ ) ~=~ \\ \\
       ~~~~~~~~~~~~=~ T^{\#} ( \inf_{X^{\#}_T}~.~.~.~ )
                                   ~=~ \inf_{Y^{\#}}~.~.~.~ ~=~ T^{\#} ( A ) ~=~ F $ \\

Thus in particular we obtain (7.5). \\

The injectivity of $T^{\#}$ will give (7.6), while (7.7) follows from (7.3) and the fact
that we can take $V = A$. \\

Conversely, let us assume (7.5). Then (7.1) collapses to the three equalities \\

$~~~ \sup_{Y^{\#}}~.~.~.~ ~=~ T^{\#} ( \sup_{X^{\#}_T}~.~.~.~ ) ~=~
            T^{\#} ( \inf_{X^{\#}_T}~.~.~.~ ) ~=~ \inf_{Y^{\#}}~.~.~.~ $ \\

thus in view of the corresponding collapsed version of (7.2), we can extend the above three
equalities to the following four \\

$~~~ \sup_{Y^{\#}}~.~.~.~ ~=~ T^{\#} ( \sup_{X^{\#}_T}~.~.~.~ ) ~=~
            T^{\#} ( \inf_{X^{\#}_T}~.~.~.~ ) ~=~ \inf_{Y^{\#}}~.~.~.~ ~=~ F $ \\

And now the injectivity of $T^{\#}$ will give (7.4) and (7.6), while (7.7) follows as
above.

\hfill $\Box$ \\

The above existence result is of a "local" nature, since it refers to a solution of the
equation (6.1) for one given right hand term $F \in Y^{\#}$. This result, however, can further
be strengthened by the following "global" one which characterizes the solvability of (6.1) for
{\it all} right hand terms $F \in Y^{\#}$. Namely, we have, [3, pp. 190,191] \\ \\

{\bf Theorem 7.2.} \\

The following are equivalent \\

(7.8)~~~ $ T^{\#} ( X^{\#}_T ) ~\supseteq~ Y $ \\

and \\

(7.9)~~~ $ T^{\#} ( X^{\#}_T ) ~=~ Y^{\#} $ \\

In each of these cases $T^{\#}$ is an order isomorphism, or in short, OI, between $X^{\#}_T$
and $Y^{\#}$. \\

It is {\it important} to note the following four facts :

\begin{itemize}

\item "Pull-back" type structures are customary when solving PDEs by functional analytic
methods. Details in this regard are presented in [3, chap. 12], while one well known classical
example has been presented in section 4 above.

\item As shown in [3, chap. 13], and indicated in section 8 next, one can consider in (5.7)
far more general partial orders than the "pull-back" type ones, and thus still obtain
solutions by the order completion method for nonlinear systems of PDEs of type (1.1), together
with possibly associated initial and/or boundary value problems.

\item Nonlinear systems of PDEs with equations of type (1.1), together with possibly
associated initial and/or boundary value problems can, in view of section 2 above, be dealt
with as equations (3.1), (3.2), and therefore (6.1). Consequently, the {\it existence} of
solutions of such systems of PDEs follows from Theorems 7.1 and 7.2.

\item The {\it blanket, type independent}, or {\it universal regularity} property of such
solutions, namely, that they can be associated with usual {\it measurable}, or even {\it
Hausdorff continuous} functions follows from arguments in [3, 1, 4-6].

\end{itemize}

{~} \\

{\bf 8. Beyond "Pull-Back" Partial Orders} \\

As indicated here, and presented in full detail in [3, chap. 13], the order completion method
in solving large classes of nonlinear systems of PDEs of the type (1.1) is {\it not} limited
to the use of "pull-back" type partial orders in (45.7), (2.3) and (2.8). In fact, a large
class of more general partial orders can be defined on the domains ${\cal M}^m_T ( \Omega )$
of the respective PDEs, and stil obtain for them solutions in the corresponding order
completions. \\

As a brief indication in this regard, let us extend the approach in section 5 above as follows.
Given any equation (3.1), (3.2) \\

(8.1)~~~ $ T ( A ) ~=~ F $ \\

where \\

(8.2)~~~ $ T : X ~\longrightarrow~ Y $ \\

is any mapping, $X$ is any nonvoid set, while $( Y, \leq )$ is a poset, while $F \in Y$ is
given, and $A \in X$ is the sought after solution. \\

Further, let $Z$ be any nonvoid set together with a {\it surjective} mapping, see (5.2) \\

(8.3)~~~ $ \lambda : Z ~\longrightarrow~ X_T $ \\

Then we obtain the commutative diagram \\

\begin{math}
\setlength{\unitlength}{0.2cm}
\thicklines
\begin{picture}(60,20)

\put(11,16){$Z$}
\put(29,18){$T_\lambda$}
\put(16,16.5){\vector(1,0){28.5}}
\put(47,16){$Y$}
\put(0,7){$(8.4)$}
\put(14,14){\vector(1,-1){12}}
\put(16,7){$\lambda$}
\put(28,0){$X_T$}
\put(33,2){\vector(1,1){12}}
\put(41,6.5){$T_\approx$}

\end{picture}
\end{math} \\

where $T_\approx$ is given in (5.3), (5.4), while we have taken by definition \\

(8.5)~~~ $ T_\lambda ~=~ T_\approx \circ \lambda $ \\

And now, we can define define the partial order $\leq_{T,~ \lambda}$ on $Z$ by \\

(8.6)~~~ $ z ~\leq_{T,~ \lambda}~ z\,' ~~~\Longleftrightarrow~~~
                      \begin{array}{|l}
                      ~~z ~=~ z\,', ~~\mbox{or} \\ \\
                      ~~ T_\lambda ( z ) \lvertneqq T_\lambda ( z\,' )
                      \end{array} $ \\

In this case, we obtain \\

(8.6)~~~ $ T_\approx ~~\mbox{is an OIE},~~~ T_\lambda ~~\mbox{is increasing},~~~
                                    \lambda ~~\mbox{is increasing and surjective} $ \\

and the procedures in sections 5 - 7 above can be applied to the mapping \\

(8.7)~~~ $ T_\lambda : Z ~~\longrightarrow~~ Y $ \\

instead of the initial mapping (3.2). \\ \\

{\bf Appendix} \\

We shortly present several notions and results used above. A related full presentation can be
found in [3, Appendix, pp. 391-420]. \\

Let $(X,\leq)$ be a nonvoid poset without minimum or maximum. For $a \in X$ we denote \\

(A.1) $~~~ < a ] = \{ x \in X ~|~ x \leq a \},~~~ [ a > = \{ x \in X ~|~ x \geq a \} $ \\

We define the mappings \\

(A.2) $~~~ X ~\supseteq~ A \longmapsto A^u = \bigcap_{a \in A}~ [ a > ~\subseteq~ X $ \\

(A.3) $~~~ X ~\supseteq~ A \longmapsto A^l = \bigcap_{a \in A} < a ] ~\subseteq~ X $ \\

then for $A \subseteq X$ we have \\

(A.4) $~~~ A^u = X \Longleftrightarrow A^l = X \Longleftrightarrow A = \phi $ \\

(A.5) $~~~ A^u = \phi \Longleftrightarrow A ~\mbox{unbounded from above} $ \\

(A.6) $~~~ A^l = \phi \Longleftrightarrow A ~\mbox{unbounded from below} $ \\

{\bf Definition A.1.} \\

We call $A \subseteq X$ a {\it cut}, if and only if \\

(A.7) $~~~A^{ul} = A $ \\

and denote \\

(A.8) $~~~ X^{\#} = \{ A \subseteq X ~|~ A ~\mbox{is a cut} \} \subseteq {\cal P} ( X) $ \\

\hfill $\Box$ \\

Clearly, (A.4) - (A.6) imply \\

(A.9) $~~~ \phi,~ X \in X^{\#} $ \\

therefore \\

(A.10) $~~~ X^{\#} \neq \phi $ \\

Given $A, B \subseteq X$, we have \\

(A.11) $~~~ A \subseteq B \Longrightarrow A^u \supseteq B^u,~ A^l \supseteq B^l $ \\

(A.12) $~~~ A \subseteq A^{ul},~~~ A \subseteq A^{lu} $ \\

(A.13) $~~~ A^{ulu} = A^u,~~~ A^{lul} = A^l $ \\

Consequently \\

(A.14) $ \begin{array}{l}
            \forall~~ A \subseteq X ~: \\ \\
            ~~~~*)~~ A^{ul} \in X^{\#} \\ \\
            \begin{array}{l}
                    ~**)~~   \forall~~ B \in X^{\#} ~: \\ \\
                       ~~~~~~~~~~~~ A \subseteq B \Longrightarrow A^{ul} \subseteq B \\ \\
                       ~~~~~~~~~~~~ B \subseteq A \Longrightarrow B \subseteq A^{ul}
                     \end{array}
          \end{array} $ \\

therefore \\

(A.15) $~~~ X^{\#} = \{ A^{ul} ~|~ A \subseteq X \} $ \\

Given $x \in X$, we have \\

(A.16) $~~~ \{ x \}^u = [ x >,~~~ \{ x \}^l = < x ],~~~ [ x >^l = < x ],~~~
                                                < x ]^u = [ x > $ \\

(A.17) $~~~ \{ x \}^{ul} = < x],~~~ \{ x \}^{lu} = [ x > $ \\

We denote for short \\

$ \{ x \}^u = x^u,~~ \{ x \}^l = x^l,~~ \{ x \}^{ul} = x^{ul},~~
                                         \{ x \}^{lu} = x^{lu},~.~.~.~ $ \\

Given $A \in X^{\#}$, we have \\

(A.18) $~~~ \phi \neq A \neq X \Longleftrightarrow
                   \left (~ \begin{array}{l}
                               \exists~~ a, b \in X ~: \\ \\
                               ~~~ < a ] ~\subseteq~ A ~\subseteq~ < b ]
                            \end{array} ~~\right ) $ \\

We shall use the {\it embedding} \\

(A.19) $~~~ X \ni x ~\stackrel{\varphi}\longmapsto~ x^{ul} = x^l = < x ] \in X^{\#} $ \\

We define on $X^{\#}$ the partial order \\

(A.20) $~~~ A \leq B \Longleftrightarrow A \subseteq B $ \\ \\

{\bf Definition 2.1.} \\

Given two posets $( X, \leq ),~ ( Y, \leq )$ and a mapping $\varphi : X \longrightarrow Y$. We
call $\varphi$ an {\it order isomorphic embedding}, or in short, OIE, if and only if it is
injective, and furthermore, for $a, b \in X$ we have \\

$~~~ a ~\leq~ b ~~~\Longleftrightarrow~~~ \varphi ( a ) ~\leq~ \varphi ( b ) $ \\

An OIE $\varphi$ is an {\it order isomorphism}, or in short, OI, if and only if it is
bijective.

\hfill $\Box$ \\

The main result concerning order completion is given in, [2] : \\ \\

{\bf Theorem ( H M MacNeille, 1937 )} \\

1)~ The poset $( X^{\#}, \leq )$ is order complete. \\

2)~ The embedding $X \stackrel{\varphi}\longrightarrow X^{\#}$ in (A.19) preserves infima and
suprema, and it is an order isomorphic embedding, or OIE. \\

3)~ For $A \in X^{\#}$, we have the order density property of $X$ in $X^{\#}$, namely \\

(A.21) $~~~ \begin{array}{l}
                A ~=~ \sup_{X^{\#}}~ \{ x^l ~|~ x \in X,~~ x^l \subseteq A \} ~=~ \\ \\
                  ~~~=~ \inf_{X^{\#}}~ \{ x^l ~|~ x \in X,~~ A \subseteq x^l \}
             \end{array} $ \\

\hfill $\Box$ \\

For $A \subseteq X$, we have \\

(A.22) $~~~ A^{ul} = \sup_{X^{\#}}~ \{ x^l ~|~ x \in A \} $ \\

Given $A_i \in X^{\#}$, with $i \in I$, we have with the partial order in $X^{\#}$ the
relations \\

(A.23) $~~~ \sup_{i \in I}~ A_i ~=~ \inf~ \{ A \in X^{\#} ~|~ \bigcup_{i \in I} A_i
                                     \subseteq A \} ~=~ ( \bigcup_{i \in I} A_i )^{ul} $ \\

(A.24) $~~~ \begin{array}{l}
                \inf_{i \in I}~ A_i ~=~ \sup~ \{ A \in X^{\#} ~|~ A \subseteq
                  \bigcap_{i \in I} A_i \} ~=~ ( \bigcap_{i \in I} A_i )^{ul} ~=~  \\ \\
               ~~~~~~~~~~~~~~=~ \bigcap_{i \in I} A_i
             \end{array} $ \\ \\

{\bf Extending mappings to order completions} \\

Let $(X,\leq),~(Y,\leq)$ be two posets without minimum or maximum, and let \\

(A.25) $~~~ \varphi : X \longrightarrow Y $ \\

be any mapping. Our interest is to obtain an extension \\

$~~~~~~ \varphi^{\#} : X^{\#} \longrightarrow Y^{\#} $ \\

For that, we first extend $\varphi$ to a {\it larger} domain, as follows \\

(A.26) $~~~ \varphi^{\#} : {\cal P} ( X ) \longrightarrow Y^{\#} $ \\

where for $A \subseteq X$ we define \\

(A.27) $~~~ \varphi^{\#} ( A ) ~=~ ( \varphi ( A ) )^{ul} ~=~
                \sup_{Y^{\#}}~ \{ < \varphi ( x )\, ] ~|~ x \in A \} $ \\

and for any mapping in (A.25), we obtain the commutative diagram \\

\begin{math}
\setlength{\unitlength}{0.1cm}
\thicklines
\begin{picture}(50,31)

\put(0,12){$(A.28)$}
\put(15,25){$X \ni x$}
\put(50,28){$\varphi$}
\put(28,26){\vector(1,0){48}}
\put(28,25){\line(0,1){2}}
\put(80,25){$\varphi ( x ) \in Y$}
\put(24.5,22){\vector(0,-1){17}}
\put(23.5,22){\line(1,0){2}}
\put(6,0){${\cal P} ( X ) \ni \{ x \}$}
\put(29,1){\vector(1,0){47}}
\put(29,0){\line(0,1){2}}
\put(80,0){$\varphi^{\#}( x ) \,=~ < \varphi ( x )\, ] \in Y^{\#}$}
\put(83,22){\vector(0,-1){17}}
\put(82,22){\line(1,0){2}}
\put(50,-3.5){$\varphi^{\#}$}

\end{picture}
\end{math} \\ \\

{\bf Proposition A.1.} \\

1)~ The mapping $\varphi^{\#} : {\cal P} ( X ) \longrightarrow Y^{\#}$ in (A.36) is
increasing, if on ${\cal P} ( X )$ we take the partial order defined by the usual inclusion
"$\subseteq$". \\

2)~ If the mapping $\varphi : X \longrightarrow Y$ in (A.35) is increasing, then the mapping
$\varphi^{\#} : {\cal P} ( X ) \longrightarrow Y^{\#}$ in (A.36) is an extension of it to
$X^{\#}$, namely, we have the commutative diagram \\

\begin{math}
\setlength{\unitlength}{0.1cm}
\thicklines
\begin{picture}(50,31)

\put(0,12){$(A.29)$}
\put(15,25){$X \ni x$}
\put(50,28){$\varphi$}
\put(28,26){\vector(1,0){48}}
\put(28,25){\line(0,1){2}}
\put(80,25){$\varphi ( x ) \in Y$}
\put(24.5,22){\vector(0,-1){17}}
\put(23.7,22){\line(1,0){2}}
\put(8,0){$X^{\#} \ni\, < x ]$}
\put(29,1){\vector(1,0){44}}
\put(29,0){\line(0,1){2}}
\put(77,0){$\varphi^{\#}( < x \,]\, ) \,=~ < \varphi ( x )\, ] \in Y^{\#}$}
\put(83,22){\vector(0,-1){17}}
\put(82,22){\line(1,0){2}}
\put(50,-4){$\varphi^{\#}$}

\end{picture}
\end{math} \\ \\

3)~ If the mapping $\varphi : X \longrightarrow Y$ in (A.25) is an OIE, then the mapping
$\varphi^{\#} : {\cal P} ( X ) \longrightarrow Y^{\#}$ in (A.26) when restricted to
$X^{\#}$, that is \\

(A.30) $~~~ \varphi^{\#} : X^{\#} \longrightarrow Y^{\#} $ \\

as in (A.29), is also an OIE. \\ \\

{\bf Lemma A.1.} \\

Let in general $\mu : M \longrightarrow N$ be an increasing mapping between two order
complete posets, then for nonvoid $E \subseteq M$ we have \\

(A.31) $~~~\mu (\, \inf_M\, E \,) ~\leq~ \inf_N\, \mu ( E )~\leq~
\sup_N\, \mu ( E )
                                                     ~\leq~ \mu (\, \sup_M\, E \,) $ \\

{\bf Proof} \\

Indeed, let $a = \inf_M\, E \in M$. Then $a \leq b$, with $b \in E$. Hence $\mu ( a ) \leq
\mu ( b )$, with $b \in E$. Thus $\mu ( a ) \leq \inf_N\, \mu ( E )$, and the first inequality
is proved. \\
The last inequality is obtained in a similar manner, while the middle inequality is trivial.

\hfill $\Box$ \\

\end{document}